\documentclass[11pt]{article}
\usepackage{latexsym}
\usepackage{amsfonts}
\makeatletter

\def\@footnotetext#1{\insert\footins{%
\footnotesize
    \interlinepenalty\interfootnotelinepenalty
    \splittopskip\footnotesep
    \splitmaxdepth \dp\strutbox \floatingpenalty \@MM
    \hsize\columnwidth \@parboxrestore
   \edef\@currentlabel{\csname p@footnote\endcsname\@thefnmark}\@makefntext
    {\rule{\z@}{\footnotesep}\ignorespaces
      #1\strut}}}
\def\abstract{\small\quotation{\hskip-\parindent\sc Abstract.}}

\def\classification{\@ifnextchar [{\@xfootnotenext}%
   {\begingroup\let\protect\noexpand
      \xdef\@thefnmark{}\endgroup
    \@footnotetext}}

\title {}

\begin{document}
\classification {{\it 1991 Mathematics
Subject Classification:} Primary 20E05,  secondary 20F28,
16S34.}

\begin{center}
{\bf \large GENERALIZED  PRIMITIVE  ELEMENTS }

{\bf \large  OF  A  FREE  GROUP}         
\bigskip

{\bf  
  Vladimir Shpilrain}

\end{center} 
\medskip

\begin{abstract}
\noindent We study endomorphisms of a free group  of finite rank 
by means of their action on specific sets of elements. In particular,
we prove that every endomorphism of the free group of rank 2 which
preserves an automorphic orbit (i.e., acts ``like  an  automorphism" 
on  one  particular orbit), is  itself an automorphism. Then, we
consider elements of a different nature, 
 defined by means of homological properties of the corresponding 
 one-relator group. 
 These elements (``generalized primitive elements"), interesting in
their own right, can  also be used for  distinguishing automorphisms 
among  arbitrary  endomorphisms. 

\end{abstract}

\date{}

\bigskip

\noindent {\bf 1.Introduction }
\bigskip

 \indent     Let  $F = F_n$   be the free group of a finite rank  
$n \ge 2$  with
a  set  
$X = {\{}x_i {\}}, 1 \le i \le n$,  of free generators.  An
element  $g \in F $ is  called {\it primitive} if it is a member of
some free basis of $ F$.  Or, equivalently:  there is an automorphism
$\phi \in Aut F$  that takes  $g$  to $ x_1 $. 

 \indent     We start here by recalling 
Problem 2 from  \cite {Sh1}.  It 
is  clear  that  the  group  $ F$   is  a  disjoint union of
different orbits under the action of the group $ Aut F$.   Denote an
orbit  ${\{}\phi(u)$,  $\phi \in Aut F{\}}$  by  $O_{Aut F}(u)$.  
It seems plausible  to  assume (see \cite {Sh1} for a general
set-up) that if an endomorphism  preserves  some orbit $ O_{Aut
F}(u)$,  i.e., if it  acts  ``like  an  automorphism"  on  one 
particular orbit, then it acts like an automorphism everywhere,
i.e., is  itself an automorphism. 

 \indent  This general conjecture appears to be quite difficult to
prove; even the case of  the  ``simplest" 
orbit  (the  one  consisting of primitive  elements)  is  still
unsettled.  Here we are able to settle the conjecture 
 in the case when $F$ has
rank 2: 
\medskip

\noindent {\bf Theorem 1.1.} If an endomorphism $\phi$ of the group
$F_2$ preserves  an  orbit $ O_{Aut F_2}(u)$, then $\phi$ is
actually an automorphism.

\medskip
 \indent Our proof of this theorem uses the following interesting
fact: there are elements of the group $F_2$ that cannot be subwords
of any {\it cyclically reduced} primitive element of $F_2$; we call
those elements {\it primitivity-blocking words}. More formally: 
\smallskip 

\noindent {\bf Definition.} An element $g \in F$ is called  a 
primitivity-blocking word if there is no cyclically reduced 
primitive element $w \in F$ such 
  that $w = g h$ for some $h \in F$ (assuming, of
course, there is no cancellation between $g$
 and $h$). 

\smallskip

\indent In the group $F_2$, primitivity-blocking words are easy
to find: for example, $[x_1, x_2]$  and $x_1^k x_2^l$ with $k, l \ge
2$, are primitivity-blocking words - this follows from a result of 
\cite {CMZ} (see also \cite {OZ}). On the other hand, in a free
group of bigger rank, the situation becomes really intriguing. We
ask: 
\medskip

\noindent {\bf Problem.} Are there primitivity-blocking words in
a free group $F_n$ if $n \ge 3$ ? 
\medskip

 \indent  
   The attempts to extend the result of Theorem 1.1 have  
 led me  to  the  following 
generalization of the primitivity concept: 
\medskip

\noindent {\bf Definition.} Let  $J$  be a right ideal of the free 
group  ring   $ZF$.   An 
element  $u \in F$  is called  $J-primitive$  if Fox derivatives \cite
{Fox}  of  the  element  $u$  generate  $J$  as a right ideal of 
$ZF$.  \smallskip

 \indent     Now ``usual" primitive elements are the same as  
$ZF$-primitive ones - this  non-trivial  fact  follows  from  a 
matrix  characterization   of  primitive elements given by
Umirbaev   \cite {U1} (cf. also \cite {Dicks}, Corollary IV.5.3). 
\medskip

 \indent     Another motivation  for  this  definition  comes  from 
commutative algebra.  We recall some facts here very briefly.  A vector  
$(p_1,...,p_m)$  
of (Laurent) polynomials is called {\it unimodular} if
$p_1,...,p_m$    generate  the whole (Laurent) polynomial algebra 
$P$  as an ideal.  Then, a  vector   $(p_1,...,p_m)$ of Laurent
polynomials is called  $\Delta-modular$  (see  \cite {A} ) 
 if  $p_1,...,p_m$ generate the augmentation ideal  $\Delta$  as an ideal of the  whole
algebra.  Suslin \cite {Su}  and Artamonov  \cite {A}   have 
proved  that  the  group   $GL_m(P), ~m \ge 3,$   acts  transitively 
on  the  set  of  all  unimodular  and  all   $\Delta$-modular 
vectors,  respectively.    See  \cite {GGN}   for   new  
interesting  applications of these results to the study of
metabelian groups.  \smallskip

 \indent        Now the above cited result of Umirbaev can  be  considered  a  free 
group analog of  Suslin's result: it says that the group  $ Aut
F$  acts  transitively on the set of all $ZF$-primitive elements
(i.e., those    with  unimodular vector of Fox derivatives). 

\smallskip

 \indent        The desire to get a free group analog of 
Artamonov's  result  draws  our attention to  $\Delta$-primitive
elements. First we note that  one can give an equivalent definition 
of $\Delta$-primitivity 
using a different  language: an element  $u \in F$ is
$\Delta$-primitive if and only if the cohomology group  $H^2(G, ZF)$
is infinite cyclic, where $G$ is one-relator group $<F ~|~ u>$. 
\smallskip

\indent 
If the group $ F$ 
has even rank  $n = 2m$, then the element   $u_n =  [x_1, x_2 ] [x_3,
x_4] ... [x_{n-1}, x_n]$   is   $\Delta$-primitive  (and  hence  so 
is  every  element from $  O_{Aut F}(u_n)$) - see Example 3.1.   On
the other hand, free group of an odd rank has no $\Delta$-primitive
elements whatsoever - see Corollary 1.5 below. R.Bieri has pointed
out to me that a combination of results of \cite{BE}, \cite{EM} and
\cite{Z} yields the following 
\medskip

\noindent {\bf Theorem 1.2.} In a group $F_{2m}$, $m \ge
1$, any $\Delta$-primitive element is   an
automorphic  image of $u_{2m}  =  [x_1, x_2 ] [x_3, x_4] ...
[x_{2m-1}, x_{2m}]$.
\medskip

\indent This, together with Corollary 1.5 below, gives a free group 
analog of Artamonov's result: 
\smallskip

\noindent {\it - in a free group $F$ of finite rank, the group $Aut
F$ acts  transitively on the set of all $\Delta$-primitive 
elements.} 

\smallskip

\indent   In the case of $F_2$, Theorem 1.2 yields   an
explicit  combinatorial
  description of $\Delta$-primitive  elements 
since every automorphic image of $[x_1, x_2 ]$ in the group $F_2$ 
has the form 
 $ [x_1,x_2 ]^g$   or $ [x_2,x_1 ]^g$   for some  $g \in F_2$.
 Note that  there  is  also 
 an explicit   combinatorial
 description of  $ZF_2$-primitive elements due to \cite {CMZ} (see
also  \cite {OZ}). 

\medskip

  \indent  There is a nice matrix  characterization of
 $\Delta$-primitive  elements; based on that, we can use 
 $\Delta$-primitive  elements for recognizing automorphisms: 
\smallskip

\noindent {\bf Proposition 1.3.}

\noindent {\bf (i)}  An element $ u \in F'_n$   is  
$\Delta$-primitive if and only if  the 
matrix $ D_u  = (d_j'(d_i (u)))_{1 \le i,j \le n} $  is
invertible over the ring $ ZF_n$ . 
\smallskip

\noindent {\bf (ii)} If $\phi$ takes a $\Delta$-primitive  
element of the group $F_n$ to another $\Delta$-primitive  element,
then $\phi$ is
 an automorphism. 

\medskip

 \indent      The matrix $ D_u$   (the ``double Jacobian" matrix) has 
been introduced  in \cite {Sh2}.   Here 
 $ d_j' $ denotes right Fox derivation whereas $ d_i $  is  the 
 ``usual",  left  Fox  derivation - see Section 2.

\medskip 
  \indent  We also treat here  the case of a free metabelian
group $M_n = F_n/F_n''$. (In the definition of a $\Delta$-primitive  
element in $M_n$, we consider abelianized Fox derivatives): 

\medskip

\noindent {\bf  Theorem 1.4.}

\noindent {\bf (i)} The group $ Aut M_2$    acts
transitively  on  the  set  of  all  $\Delta$-primitive elements of
$ M_2$.  In other words, every  $\Delta$-primitive element of
$M_2$   is an automorphic image of the element  $[x _1,x_2 ]$  and 
therefore  has the form $ [x _1,x_2 ]^g$   or $ [x_2,x_1 ]^g$   for
some  $g \in M_2$ ; \smallskip

\noindent {\bf (ii)} There are no $\Delta$-primitive elements in
the group $M_n$ if $n$ is odd.

\smallskip 

\indent Part (i) of this  theorem answers a question  from \cite
{GGN} in rank 2 case.  
\medskip

\noindent {\bf Corollary 1.5.} There are no $\Delta$-primitive
elements in the free group $F_n$ if $n$ is odd. 

\medskip

\indent  Talking  about  general  case  of   $J$-primitive 
elements  for   an arbitrary right ideal  $J$ of the group ring $ZF$,   
we  note  first  of  all  that   $J$-primitive  elements do not
always exist.  For example, if  $J$  cannot  be  generated  by 
less  than   $(n + 1)$   elements,  then  there  are  obviously  
no   $J$-primitive  elements  in  the  group  $ F_n$.   In  fact, 
if  $ g $  is   a  $J$-primitive element, and the ideal  $J$  is  
$k$-generated,  then   $g$   has {\it outer rank} $ k$  in the sense
defined  in \cite {Sh1} (the  minimal  number  of  free  generators
on which an automorphic  image  of  $ g$   can  depend).   This 
follows from a result of Umirbaev \cite {U2}. 

 \indent     Deciding 
  for which (right) ideals $J$ of the group ring $ZF$ 
 the group  $Aut F$  acts transitively 
on  the  set  of all  $J$-primitive elements of $F$, seems to be a
difficult problem.
 Actually, it makes sense only for {\it
characteristic}  ideals, i.e., those invariant under free group
automorphisms. 
 Note that if an element $g \in F$ is $J$-primitive for a 
 characteristic right ideal  $J$, then every automorphic image
of $g$ is $J$-primitive, too - this follows from the equality (1)
in the next section. 
\smallskip

 \indent  Also, it is easy to show   that for any $m$-generator
right ideal  $J$  of   $ZF$,  the group $ GL_m (ZF)$  acts
transitively on the set of all $m$-tuples  of elements generating 
$J$  (Lemma 2.4).  \smallskip

 \indent         Of particular interest  are  $J$-primitive
elements  of  maximal outer rank $ n = rank F$  in the case when  
$J$  is the right ideal  generated  by Fox derivatives of a group
element: 

\smallskip
\noindent {\bf Proposition 1.6.} Let  $g \in F_n$ be an
element of outer rank $n$,  and let $J$   be the  right  ideal  of 
$ZF_n$   generated by Fox derivatives of $ g$.  If a {\it
monomorphism }  (i.e., injective     endomorphism) $\phi$  of the 
group   $F_n$  takes  $g$  to another  $J$-primitive element, then 
$\phi$  is  actually  an  automorphism. 
\smallskip

 \indent     In particular: 

\smallskip

\noindent {\bf Corollary 1.7.} (cf. \cite {T}). Let  $g \in F_n$  
be an element  of  outer rank $ n$.  If $ \phi(g) = g$  for some
monomorphism  $\phi$  of  the  group   $F_n$ ,   then   $\phi$   is
actually an automorphism. \\

\noindent {\bf 2. Preliminaries  }
\bigskip

  \indent      Let  $ZF$  be the integral group ring of the 
group  $F$ 
and  $\Delta$   its augmentation ideal, that is, the kernel of the
natural homomorphism  $\epsilon$ : 
$ZF \to Z$.  More generally, when  $R \subseteq F$  is a normal 
subgroup  of   $F$,  we denote by $\Delta_R$   the ideal of  $ZF$ 
generated by all elements  of  the  form  
$(r - 1)$,  $r \in R$.  It is the kernel of the natural homomorphism 
$\epsilon_R$ : $ZF \to Z(F/R)$. 
\smallskip

  \indent     The ideal  $\Delta$  is a free left  $ZF$-module with
a free basis ${\{}(x_i  - 1){\}}$,  $1\le i \le n$, and left Fox
derivations $ d_i $  are projections  to 
the corresponding free cyclic direct summands.  Thus any element 
$ u \in \Delta$ can be uniquely written in the form  
$u =  \sum_{i=1}^n  d_i (u) (x_i - 1)$. 

 \indent     Since the ideal  $\Delta$  is a free right 
$ZF$-module  as  well,  one  can define right Fox derivatives 
$d_i'(u)$  accordingly, so that  $u =  \sum_{i=1}^n  (x_i  - 1)
d_i'(u)$. 

 \indent        One  can  extend  these  derivations  linearly  to 
 the  whole   $ZF$  by setting $ d_i'(1) = d_i (1) = 0$.  
\smallskip

 \indent       The next lemma is an immediate consequence  of the
 definitions.

\smallskip

\noindent {\bf Lemma 2.1.}  Let  $J$  be an arbitrary left (right)
ideal of  $ZF$   and  let  $u \in \Delta$.  Then $u \in \Delta J$ 
 ($u \in J\Delta$ )  if and only if $d_i'(u) \in J $ ($d_i(u) \in
J$)  for each  $i, 1 \le i \le n$. 

\smallskip

 \indent      Proof of the next lemma can be found in \cite {Fox}. 
\smallskip

\noindent {\bf Lemma 2.2.}  Let  $R$  be a normal subgroup of  $F$,
and let  $y \in F$.  Then   $y - 1 \in$ 

\noindent 
 $\Delta_R \Delta$  if and
only if  $y \in R'$. 
\smallskip

 \indent      We also need the ``chain rule" for Fox derivations
(see \cite {Fox}): 
\smallskip

\noindent {\bf Lemma 2.3.} Let  $\phi$  be  an  endomorphism  of  
$F$   (it  can  be  linearly extended to  $ZF$) defined by 
$\phi(x_k) = y_k , 1 \le k \le n$,  and let  $v =  \phi(u) $ 
for some  $u, v \in ZF$.  Then:
\begin{center}

     $ d_j(v) =   \sum_{k=1}^n \phi(d_k(u))d_j(y_k)$.   
\end{center} 
                   
\smallskip

\indent      For an endomorphism  $\phi : x_i \to y_i $ , $1 \le
     i \le n$ ,   of  the  group    $F_n$, 
let  $J_{\phi}  = (d_j (y_i))_{1 \le i,j \le n}$     be the
Jacobian matrix of $\phi$.  We are going to 
need the following application of Lemma 2.3: 
 if  $g, h \in  F_n$   and   $h  = \phi(g)$,  then 

\begin{eqnarray}
(d_1(h)),...,d_n(h))=(\phi(d_1(g)),...,\phi(d_n(g)))J_\phi.         
\end{eqnarray} 

\indent      Now comes the key lemma: 

\smallskip

\noindent {\bf Lemma 2.4.} Let  $J$  be a right ideal of  $ZF$ 
(hence a free  right  module over  $ZF$ - see \cite {C}) generated
as a free module over  $ZF$  by   $u_1 ,...,u_m $.  
Then the following conditions are equivalent: 
\smallskip

\noindent {\bf (a)} A matrix $M = \left( a_{ij}\right)_{1\le i, j\le
m}$   is invertible
over  $ZF$  (i.e.,  $M  \in GL_m(ZF)$) ; 

\medskip

\noindent {\bf (b)} The elements  $y_j = \sum_{k=1}^m  (u_k - 1)
a_{kj}$, $1\le j\le m$,
 generate the ideal  $J$  as a right ideal of  $ZF$. 
\smallskip

\noindent {\bf Proof.} Suppose  $M$  is invertible over $ ZF$ ;
  denote by  $U$  the row matrix  
$(u_1 ,...,u_m)$,  and by  $Y$ - the row matrix  $(y_1,...,y_m)$. 
Then  $U M M^{-1}   = Y M^{-1}   = U$  which means that $u_1 ,...,u_m$ 
belong to the right ideal of  $ZF $ generated by $ y_1,...,y_m$ . 
Conversely, suppose we have $ Y B = U $ for some matrix  $B $ over 
$ZF$.  Then  $U M B = U$,  hence  $M B  =  I$,  the  identity
matrix, because    $(u_1 ,...,u_m) $  form  a  free  basis  of  a 
free  right  $ZF$-module  $J$.  This implies  $M \in GL_m(ZF)$ - see
e.g. \cite {C}. \\

\noindent {\bf 3. Proofs}
\bigskip

\noindent {\bf Proof of Theorem 1.1.} We consider several
possibilities for an orbit $ O_{Aut F_2}(u)$: 
\smallskip

\noindent (1) $u$ is a primitive element. In this case, the image
of every primitive element of $F_2$ under the  endomorphism $\phi$
is primitive. By composing $\phi$ with some automorphism of the
group $F_2$ if necessary, we may assume that $\phi(x_1) = x_1$. 

\indent Now write $\phi(x_2) = x_1^k g$, where $k$ is an integer, 
and $g$ is an element of the normal closure of $x_2$, so that $g$
is a product of elements of the form $h_i  x_2^{\pm1} h_i^{-1}$.

\indent Since every element of the form $x_1^{n} x_2$  is primitive, 
its image $s = \phi(x_1^{n} x_2) = x_1^{n+k} g$ is
primitive, too. 

\indent Recall now a result of \cite{CMZ} which says that if $w$ is
a primitive element of $F_2$, then some conjugate of $w$ can be
written in the form $x_1^{k_1} x_2^{l_1} ... x_1^{k_m} x_2^{l_m}$,
so that some of $x_i$ occurs either solely in the exponent 1 or 
solely in the exponent $-$1. 

\indent 
We see that for sufficiently large $n$, 
generator $x_1$ would not occur in any conjugate of 
$s$ neither solely in the exponent 1 nor solely in the exponent $-$1.
Therefore, $x_2$ should be the one with this property. 

\indent It follows that  no $h_i$ in the  decomposition of $g$
mentioned above, has  entries of $x_2$, so that every $h_i$ is a 
power of $x_1$. 

\indent Then, since 
the element $x_1^{-k }x_2 $ is primitive, its image
$\phi(x_1^{-k}x_2) = g$ is primitive as well. This implies that the
sum of powers of $x_2$ in the  decomposition of $g$ should be equal
to $\pm 1$; otherwise $g$ would not be primitive even modulo
$F_2'$.  It follows that our $g$ is actually 
conjugate to $x_2^{\pm1}$. 

\indent Summing up, we see that $g$ has the form $h x_2^{\pm1}
h^{-1}$ for  some $h \in F_2$, $h$ a power of $x_1$.  
Therefore,  $\phi(x_1) =  x_1$ 
  and $\phi(x_2) = x_1^k g$ generate  the group 
$F_2$,  i.e., $\phi$ is an automorphism.
\smallskip

\noindent (2) $u = v^k$ is a power of a primitive element $v$. In
this case, every image
of  $u $ under the  endomorphism $\phi$ has the form $w^k$, $w$ a
primitive element. It follows that the image
of every primitive element of $F_2$ is primitive, hence we may
apply the argument from the previous case. 
 
\smallskip

\noindent (3) $u$ has outer rank 2, i.e, $u$ does not belong to a
proper free factor of $F_2$. 
If $u$ does not belong to a proper $retract$ of $F_2$, then $\phi$ 
 is an automorphism by a result of Turner \cite {T}. 

\indent Suppose now that $u$ belongs to a  proper retract $R$ of
$F_2$, and $\phi$ is the corresponding $retraction$, i.e., 
$\phi(F_2) = R$ (otherwise $\phi$ would already be an automorphism 
by \cite {T}). Since $R$ is a  proper retract of $F_2$, it should 
have rank smaller than 2, i.e., $R$ is cyclic,  and $R \subseteq 
O_{Aut F}(u)$. 
\smallskip

\indent We show now that no automorphic orbit $O_{Aut F}(u)$ can 
contain a non-trivial cyclic group. By means of contradiction,
suppose some automorphic image of some $s \in O_{Aut F}(u)$  is of
the form $s^k$, $k \ge 2$,
 and let $s$ have minimal length among all the
elements of $O_{Aut F}(u)$ with this property. Let $\alpha(s) = s^k$ 
with $k \ge 2$, and let  $\alpha^{-1}(s) = r$. Then $\alpha(s) =
\alpha(r)^k = \alpha(r^k)$,  hence $s = r^k$. Thus, every
automorphic image of  $r$ is of the form $r^m$; furthermore, $r \in
O_{Aut F}(u)$ (since $r = \alpha^{-1}(s)$), and $r$ has length
smaller than  that of $s$ (see \cite{LS}, Proposition I.2.15).
This contradiction completes the proof of Theorem 1.1.

\medskip
\indent Before getting to a proof of Theorem 1.2, we consider a
couple of examples. 
\medskip

\noindent {\bf Example 3.1.} The element $u =  [x_1, x_2 ] [x_3, x_4] 
... [x_{m-1}, x_m]$ of the group $F_{2m}$  is  $\Delta$-primitive  
since the corresponding double Jacobian matrix $ D_u $ 
 is invertible - see \cite {Sh1}, Proposition 4.1. 
\smallskip

\noindent {\bf Example 3.2.} The element $v =  [x_1, x_2 ] [x_2,
x_3] [x_3, x_4]$ of the group $F_4$  is  $\Delta$-primitive. It is
not quite obvious that $v$  is an automorphic image of $u = [x _1,x_2
] [x_3, x_4]$. However, this is the case: $u$ is taken to $v$ by
the following automorphism: $x_1 \to x_1 x_3^{-1}$; $x_2 \to x_3 x_2 
x_3^{-1}$; $x_3 \to x_3$; $x_4 \to x_4$. 

\smallskip
\indent Now we get to Theorem 1.2: 

\medskip
\noindent {\bf Proof of Theorem 1.2.} We give a proof here without 
introducing background material, just referring to \cite{BE} and 
\cite{EM} for details. 

\indent First of all, it is an immediate consequence  of the
 definition that  $g \in F$ is a $\Delta$-primitive element
if and only if for the right ideal $J_g$ of the group ring $ZF$ 
generated by Fox derivatives of $u$, one has factor-module 
$ZF/J_g$ isomorphic to the trivial $F$-module $Z$. In other words,
as we have mentioned in the Introduction, the cohomology group  
$H^2(G, ZF)$ of 
  the one-relator group $G = <F ~|~ g>$  is  infinite cyclic. 
It is also clear that the group $G$ is torsion-free. 

\indent Then, Theorem 9.3 of \cite{BE} implies that the  
group $<F ~|~ g>$  is a {\it Poincar\'e 
 duality group of dimension 2} ($PD^2$-group). It has been proved 
later in \cite{EM} that every one-relator (torsion-free) 
$PD^2$-group is a surface group. 

\indent Applying now a result of Zieschang \cite{Z}, we see that,
in case $F = F_{2m}$, $g$ must be an
automorphic  image of $u_{2m}  =  [x_1, x_2 ] [x_3, x_4] ...
[x_{2m-1}, x_{2m}]$. This completes the proof. 

\medskip

\noindent {\bf Proof of Proposition 1.3.} 

\noindent (i) First suppose the matrix  $ D_u  = (d_j'(d_i(u)))_{1
\le i,j \le n} $    is invertible. Then, by Lemma 2.4, the elements 
$y_i = \sum_{k=1}^n  (x_k - 1)d_k'(d_i(u)) =  d_i(u) -
\epsilon(d_i (u))$,  $1\le i\le n$,
 generate the ideal  $\Delta$  as a right ideal of  $ZF$. Since $u
\in F'$, we have $\epsilon(d_i (u)) = 0$, $1\le i\le n$,
by Lemmas 2.1, 2.2. Therefore, $u$ is $\Delta$-primitive.

\indent Conversely, if $u$ is $\Delta$-primitive, then the elements
$d_i(u) =  \sum_{k=1}^n  (x_k - 1)d_k'(d_i(u))$, generate 
$\Delta$  as a right ideal of  $ZF$. Again  by Lemma 2.4,
the matrix  $ D_u $ is invertible.

\smallskip
\noindent (ii) We use an 
 argument from \cite{Sh2} implying
that if a matrix $ D_{\phi(u)} $ is invertible and $u \in F'$, then
$\phi$  is an automorphism. Applying part (i) of this proposition 
completes the proof. 

 \medskip

\noindent {\bf Proof of Theorem 1.4.} First of all, note that
if $u$ is $\Delta$-primitive, then 
$u \in M'$; otherwise some  $d_i(u)$ wouldn't belong to $\Delta$ 
by Lemmas 2.1, 2.2.
\smallskip

\noindent {\bf (i)} Let $h \in M_2'$; then we can write $h$ as 
$[x_1, x_2]^w $ for some $w \in  ZA_2 = Z(M_2/M_2')$. Then for
abelianized Fox derivatives (we denote them the same way as the ones
in a free group
 ring when there is no ambiguity) we have: 
$ d_i(h) = w \cdot d_i([x_1, x_2])$. Hence  if
 $d_i(h), i = 1,2,$ generate the same ideal of  $ZA_2$ as $d_i([x_1,
x_2])$, the element $w$ should be invertible in $ZA_2$, which means
it has the form $\pm g$ for some $g \in M_2/M_2'$. 

\indent Thus  $h = [x_1, x_2]^{\pm g}$; in particular, it  is an
automorphic image  of $[x_1, x_2]$. 
\smallskip

\noindent {\bf (ii)} Consider basic commutators of weight 2 in the
group $M_n$:  $c_1 = [x_1, x_2]$; $c_2 = [x_1, x_3]$; ...;  $c_N = 
[x_{n-1}, x_n]$, where $N = n(n - 1)/2$, 
 and consider a product of the form $w = c_1^{k_1} 
c_2^{k_2} ... c_N^{k_N}$.
  Evaluate   abelianized Fox derivatives of the element $w$: 

\smallskip
\noindent $d_1(w) = k_1(x_2 - 1) + ... + k_{n - 1}(x_n - 1)$;

\noindent $d_2(w) = -k_1(x_1 - 1) + k_n(x_3 - 1) + ... +  k_{2n -
3}(x_n - 1)$; 

\noindent $d_3(w) = -k_2(x_1 - 1) - k_n(x_2 - 1) +
k_{2n - 2}(x_4 - 1) + ... + k_{3n - 6}(x_n - 1)$. 

\begin{center}
. . .
\end{center} 

\noindent $d_n(w) = -k_{n - 1}(x_1 - 1) - ... - k_{n(n - 1)/2}(x_{n -
1} - 1)$. 
\smallskip

\indent We are now going to show that these derivatives do not
generate $\Delta$ even modulo $\Delta^2$. To do that, it suffices
to show that they are linearly dependent, i.e., that 
the $n$x$n$ matrix of coefficients (its $(i, j)$th entry
is the coefficient at  $(x_j - 1)$ in the decomposition of $d_i(w)$
above) has determinant 0.  

\indent It is easy to see that this matrix (denote it by $A =
(a_{ij}$)) is antisymmetric, with zeroes on the diagonal. The
determinant of a matrix like that must be 0 if $n$ is odd. Indeed,
consider a summand $a_{1,i_1} a_{2,i_2} ... a_{n,i_n}$ in the
decomposition of the determinant. If there is at least one diagonal
element among these $a_{k,i_k}$, then the product is 0. If all 
$a_{k,i_k}$ are off-diagonal elements,  consider the ``reflection" 
$a_{i_1,1} a_{i_2,2} ... a_{i_n,n}$. These 2 summands  go with
different signs since $a_{ij} = - a_{ji}$, and $n$ is odd.
Therefore, they cancel out which proves that the determinant of $A$
equals 0. This completes the proof of part (ii). 
\medskip

\noindent {\bf Proof of Corollary 1.5.} If there were a
$\Delta$-primitive element  in  $F_n$, its image in  $M_n$ would be a
$\Delta$-primitive  element  of $M_n$ which contradicts Theorem 
1.4 (ii). 

\medskip
\noindent {\bf Proof of Proposition 1.6.} Let $h  = \phi(g)$ be a
$J$-primitive element. Since $g$ has outer rank $n$, the right ideal  $J$
is $n$-generated by \cite {U2}. It follows that the elements 
$d_1(g),...,d_n(g)$ freely generate $J$ as a right ideal of $ZF_n$.
Indeed, if there were a (right) $ZF_n$-dependence between these 
elements, then one of them would belong to the right ideal
generated by the others - this follows from a general theory of
\cite {C}. Therefore, $J$ could be generated by less than $n$
elements. 

\indent  Thus by Lemma 2.4,  for some matrix  $M \in GL_n(ZF_n)$, we
have: 
\begin{eqnarray}
(d_1(h)),...,d_n(h)) = (d_1(g)),...,d_n(g)) M. 
\end{eqnarray} 

\indent  On the other hand, by the equation (1) in the Preliminaries,
we have: 

\begin{center}

$(d_1(h),...,d_n(h))$ =
$(\phi(d_1(g)),...,\phi(d_n(g))) J_\phi.$        
   
\end{center} 
\smallskip

\indent This together with (2) gives 

\begin{eqnarray}
(d_1(g),...,d_n(g)) = (\phi(d_1(g)),...,\phi(d_n(g)))J_\phi
M^{-1}.
\end{eqnarray} 

\indent  This means  $J \subseteq  \phi(J)$. Since $\phi$ is a
monomorphism, this yields   $J = \phi(J)$, in which case the matrix
$J_\phi M^{-1}$ on the right-hand side of (3) must be invertible by
Lemma 2.4.  Therefore, $J_\phi$ is invertible, too, hence $\phi \in
Aut F_n$ by \cite {B}.

\medskip
\noindent {\bf Corollary 1.7} follows immediately.\\

\noindent {\bf Acknowledgement}
\smallskip

\indent I am grateful to S.Ivanov for numerous useful discussions, 
and to R.Bieri - for pointing out Theorem 1.2.

\medskip
\noindent 
 Department of Mathematics, University of California, 
Santa Barbara, CA 93106 
 
\smallskip

\noindent {\it e-mail address\/}: shpil@math.ucsb.edu

\end{document}